\documentclass{article}                                       

\def\thebibliography#1{\small                                 
       \list{\arabic{enumi}.}                                 
        {\settowidth\labelwidth{[#1]}                         
        \leftmargin\labelwidth                                
        \itemsep 0pt                                          
        \parsep \itemsep                                      
        \advance\leftmargin\labelsep                          
        \usecounter{enumi}}}                                  
\def\@cite#1{$^{#1}$} 
\def\mytitle#1{\setcounter{equation}{0}
\setcounter{footnote}{0}
\begin{flushleft}{\bf\textbf{\Large #1}}\end{flushleft}       
\vspace{0.25cm}}                                              
\def\myname#1{\leftline{{\large #1}}\vspace{-0.13cm}}                                                              
\def\myplace#1#2{\small\begin{flushleft}{\it #1}\\
{\tt#2}\end{flushleft}}                    
\newenvironment{myabstract}{\normalsize\noindent}{}           

\def\myclassification#1{\small\noindent 2000 Mathematics Subject Classification. #1\vspace{0.5cm}} 
\begin{document}                                              

\mytitle{Matrix Problems in Hilbert Spaces}


\myname{A.V. Roiter, S.A. Kruglyak, L.A. Nazarova}

\myplace{Institute of Mathematics, 01601, Tereshchenkovska str., 3, Kiev, Ukraine}
{[roiter@imath.kiev.ua]}


\myclassification{16G20}

\begin{myabstract}
Thirty years ago (see~\cite{1}) it was shown in the works of
I.M.~Gelfand, V.A.~Ponomarev, P.~Gabriel and the authors that a number of
problems of linear algebra admits deep research on the one hand
in naive terms of reducing of collections of matrices by means of certain
admissible transformations, on the other hand in categorial terms.
Hereat the same problems have arisen in representations of algebras,
group theory and Harish-Chandra modules. Firstly, representations of
quivers and posets should be mentioned here. It was proved in particular that
finitely represented (resp. tame) quivers correspond to Dynkin (resp.
extended Dynkin) graphs. A way of transmitting of these results to Hilbert
spaces was shown in \cite{2}. Coxeter functors also were constructed and an
analogue of Gabriel theorem \cite{3} was proved there.

Put $Q_v$ be the set of vertices, $Q_a$ be the set of arrows of quiver
(= oriented graph) $Q$, $q=|Q_v|<\infty$; $t_\alpha$ is tail, $h_\alpha$
is head of $\alpha\in Q_a$. For $a\in Q_v$ put $T(a)=\{\alpha\in
Q_a\,|\,t_\alpha =a\}$, $H(a)=\{\alpha\in Q_a\,|\,h_\alpha = a\}$,
$Q(a,b)=H(b)\cap T(a)$ ($a,b\in Q_v$). Non-oriented graph $\Gamma(Q)$
naturally corresponds to quiver $Q$.

Representation $S$ of quiver $Q$ in category $K$
attaches $S(a)\in \mathrm{Ob}\, K$ for any $a\in Q_v$ and
morphism $S(\alpha): S(a)\to S(b)$ for any $\alpha\in Q(a,b)$.
Let $\mathcal{H}$ be a category, which objects are
Hilbert spaces and morphisms are bounded linear
operators. $\varphi^*\in\mathcal{H}(B,A)$ uniquely corresponds to
any $\varphi\in\mathcal{H}(A,B)$.
A real nonzero nonnegative function $\chi$ on $Q_v$ is a
\emph{character} on $Q$.
We assume here that $\Gamma(Q)$ is a tree. A representation $S$
is \emph{orthoscalar} if for some character $\chi_S$ and any $a\in Q_v$
$$
\sum_{\alpha\in
T(a)} S^*(\alpha)S(\alpha)+ \sum_{\beta\in H(a)}
S(\beta)S^*(\beta)= \chi_S(a){\bf 1}_a.
$$
It may be proved that {\it if $\Gamma(Q)$ is neither
Dynkin graph \cite{2} nor extended Dynkin graph then the problem of
classification of indecomposable orthoscalar representations of $Q$ is
wild. Let $\Gamma(Q)$ be an extended Dynkin graph. $S$ is an
orthoscalar indecomposable faithful ($S(\alpha)\neq 0$, $\alpha\in Q_a$)
representation. Then $d=\dim(S)=(\dim S(v_1),\ldots,\dim S(v_q))$ is either
a real or the minimal imaginary root of $\Gamma(Q)$ ($\dim
S(v_i)<\infty$, see \cite{4}). Given $d$ such representations are discribed
up to unitary equivalence in the first case by $q-1$ and in the second one
by $q+1$ real parameters.}

Orthoscalar representations of posets and other matrix problems also arise.
\end{myabstract}

\end{document}